\documentclass[a4paper]{amsart}

\usepackage[mathscr]{eucal}
\usepackage{amsmath,amsfonts,amssymb}
\usepackage{mathrsfs}
\usepackage[dvips]{graphicx,color}
\usepackage{fancyhdr}

\newcommand{\mathsym}[1]{{}}

\theoremstyle{plain}

\newtheorem{theorem}{Theorem}[section]

\newtheorem{proposition}{Proposition}[section]

\newtheorem{remark}{Remark}[section]

\makeatletter
\@addtoreset{equation}{section}
\makeatother

\begin{document}

\title{A note on\\ zeros of Eisenstein series for genus zero Fuchsian groups}
\author{Junichi Shigezumi}
\date{}

\maketitle\vspace{-0.2in}

\begin{abstract}
Let $\Gamma \subseteq \text{SL}_2(\mathbb{R})$ be a genus zero Fuchsian group of the first kind having $\infty$ as a cusp, and let $E_{2 k}^{\Gamma}$ be the holomorphic Eisenstein series associated with $\Gamma$ for the $\infty$ cusp that does not vanish at $\infty$ but vanishes at all the other cusps. In the paper ``On zeros of Eisenstein series for genus zero Fuchsian groups'', under assumptions on $\Gamma$, and on a certain fundamental domain $\mathcal{F}$, H. Hahn proved that all but at most $c(\Gamma, \mathcal{F})$ (a constant) of the zeros of $E_{2 k}^{\Gamma}$ lie on a certain subset of $\{ z \in \mathfrak{H} \, : \, j_{\Gamma}(z) \in \mathbb{R}\}$.

In this note, we consider a small generalization of Hahn's result on the domain locating the zeros of $E_{2 k}^{\Gamma}$. We can prove most of the zeros of $E_{2 k}^{\Gamma}$ in $\mathcal{F}$ lie on its lower arcs under the same assumption.
\end{abstract}\quad

\markboth{Note on zeros of Eisenstein series for genus zero Fuchsian groups}{Note on zeros of Eisenstein series for genus zero Fuchsian groups}

\section{Introduction}
We denote by $\Gamma$ a Fuchsian group of the first kind, which has $\infty$ as a cusp with width $h$. Let $\mathcal{F}$ be a fundamental domain of $\Gamma$ contained in $\{ z \, : \, -h/2 \leqslant \Re(z) < h/2 \}$ and $A$ be the lower arcs of $\mathcal{F}$. We then define
\begin{align}
y_0 &:= \inf\{ y \, : \, \pm h/2 + i \, y \in \partial \mathcal{F} \},\\
a_0 &:= j_{\Gamma}(-h/2 + i \, y_0).\label{a0}
\end{align}
We denote by $c(\Gamma, \mathcal{F})$ the number of equivalence classes under the action of $\Gamma$ on the set of critical points of $\mathcal{F}$ at which $\frac{y'(t)}{(j_{\Gamma} \circ z_{A})'(t)}$ changes sign. Further details may be found in \cite{H}.

The main theorem of \cite{H} is the following:

\begin{theorem}\cite[Theorem 1.1]{H}\quad
Let $\Gamma$ be a genus zero group that is good for the weight $2 k$. Suppose that $\mathcal{F}$ is acceptable for $\Gamma$ and that $j_{\Gamma}$ has real Fourier coefficients. Then all but possibly $c(\Gamma, \mathcal{F})$ of the zeros of $P(E_{2 k}^{\Gamma}, X)$ lie on $[a_0, \infty)$, where $a_0$ is as in $(\ref{a0})$. Moreover, if $m$ denotes the number of distinct zeros with odd multiplicity on $[a_0, \infty)$, then $m + c(\Gamma, \mathcal{F}) \geqslant \deg(P(E_{2 k}^{\Gamma}, X))$.\label{th-main0}
\end{theorem}
Here, {\it ``good''} and {\it ``acceptable'}' are defined in \cite[Section 1]{H}, and $P(f, X) \in \mathbb{C}[X]$ is the {\it divisor polynomial} of a modular form $f$. (see \cite[Section 3]{H})\\

Now, we denote by $s_{2 k}^1(\Gamma)$ the number of cusps other than $0$ and $-1/2$ at which $E_{2 k}^{\Gamma}$ has odd multiplicity. Then, since $E_{2 k}^{\Gamma}$ vanishes at all of the cusps other than $\infty$, we can substitute $c(\Gamma, \mathcal{F})$ for $c(\Gamma, \mathcal{F}) - s_{2 k}^1(\Gamma)$ in Theorem \ref{th-main0} (see Section \ref{sec-ex-c}). Similarly to $y_0$ and $a_0$, we define
\begin{align}
y_1 &:= \inf\{ y \, : \, i \, y \in \partial \mathcal{F} \},\label{y1}\\
a_1 &:= j_{\Gamma}(i \, y_1).\label{a1}
\end{align}

Then we state the following theorem:

\begin{theorem}
Let $\Gamma$ be a genus zero group that is good for the weight $2 k$. Suppose that $\mathcal{F}$ is acceptable for $\Gamma$ and that $j_{\Gamma}$ has real Fourier coefficients. Then all but at most $2 (c(\Gamma, \mathcal{F}) - s_{2 k}^1(\Gamma))$ of the zeros of $P(E_{2 k}^{\Gamma}, X)$ lie on $[a_0, a_1]$.\label{th-main}
\end{theorem}\quad

For the examples of \cite[Section 2]{H}, if $\Gamma = \text{SL}_2(\mathbb{Z})$, then $c(\Gamma, \mathcal{F}) - s_{2 k}^1(\Gamma) = 0$, and so we have proved that all of the zeros of $E_{2 k}^{\Gamma}$ in $\mathcal{F}$ lie on its lower arcs for this cases. On the other hand, if $\Gamma = \Gamma_0(3)$, then $c(\Gamma, \mathcal{F}) - s_{2 k}^1(\Gamma) = 1$, and so we have proved that all but at most two of the zeros of $E_{2 k}^{\Gamma}$ in $\mathcal{F}$ lie on its lower arcs.
\newpage

\section{An expansion on the number $c(\Gamma, \mathcal{F})$}\label{sec-ex-c}
Let the group $\Gamma$ has an acceptable fundamental domain. Note that, other than $\infty$, $0$, and $-1/2$, all of the cusps are critical points at which $\frac{\partial y}{\partial j_{\Gamma}}$ changes sign. Again, $E_{2 k}^{\Gamma}$ vanishes at every cusp except $\infty$. Thus, $E_{2 k}^{\Gamma}$ has $s(\Gamma) - s_{2 k}^1(\Gamma)$ zeros at cusps whose order is even and at least $2$.

Let $\mathfrak{a}_1$, $\ldots$ , $\mathfrak{a}_{m'}$ be the zeros of $P(E_{2 k}^{\Gamma}, X)$ on $[a_0, \infty)$ excluding images of cusps by $j_{\Gamma}$ that have odd multiplicity, and let $\mathfrak{b}_1$, $\ldots$ , $\mathfrak{b}_{c(\Gamma, \mathcal{F}) - s_{2 k}^1(\Gamma, \mathcal{F})}$ be the points excluding cusps along the lower arcs of $\mathcal{F}$ where $\frac{\partial y}{\partial j_{\Gamma}}$ changes sign. Let
\begin{equation}
Q_1(X) := \prod_{i = 1}^{m'} (X - \mathfrak{a}_i) \prod_{j = 1}^{c(\Gamma, \mathcal{F}) - s_{2 k}^1(\Gamma, \mathcal{F})} (X - j_{\Gamma}(\mathfrak{b}_j)).
\end{equation}

Then, similarly to the proof of Theorem \ref{th-main0} (cf. \cite[Proof of Theorem 1.1, Section 4]{H}), $Q_1(j_{\Gamma}) \Upsilon \in M_{2 k}^{\infty}(\Gamma)$ and $E_{2 k}^{\Gamma}$ satisfy
\begin{equation}
\langle E_{2 k}^{\Gamma}, Q_1(j_{\Gamma}) \Upsilon \rangle = 0,\label{ort-EQd0}
\end{equation}
and so satisfy the assumption of Proposition 4.1 of \cite{H}. Then, we assume $m + c(\Gamma, \mathcal{F}) - s_{2 k}^1(\Gamma) < \deg(P(E_{2 k}^{\Gamma}, X))$, that is $m' + c(\Gamma, \mathcal{F}) - s_{2 k}^1(\Gamma) < \deg(P(E_{2 k}^{\Gamma}, X)) - s_{2 k}^1(\Gamma) \leqslant \deg(P(E_{2 k}^{\Gamma}, X))$. We have
\begin{equation}
\Re (\langle E_{2 k}^{\Gamma}, Q_1(j_{\Gamma}) \Upsilon \rangle)
 = - \sum_{i = 1}^{n} \int_{A_i} \left( \int_{\infty}^{j_{\Gamma}(z)} P(E_{2 k}^{\Gamma}, j_{\Gamma}) \, Q_1(j_{\Gamma}) \, (|\Upsilon|^2 t^{2 k}) (j_{\Gamma}) \, t^{-2} (j_{\Gamma}) \, t' (j_{\Gamma}) \, d j_{\Gamma} \right) d x.
\end{equation}
However, this contradicts relation (\ref{ort-EQd0}) since $P(E_{2 k}^{\Gamma}, X) Q_1(X) t'(X)$ is nonnegative or nonpositive on $[j_{\Gamma}(z), \infty)$.

 We can thus state the following theorem in place of Theorem \ref{th-main0}:

\quad \vspace{-0.07in} \\
{\bfseries Theorem \ref{th-main0}.$'$} {\itshape Let $\Gamma$ be a genus zero group that is good for the weight $2 k$. Suppose that $\mathcal{F}$ is acceptable for $\Gamma$ and that $j_{\Gamma}$ has real Fourier coefficients. Then all but possibly $c(\Gamma, \mathcal{F}) - s_{2 k}^1(\Gamma)$ of the zeros of $P(E_{2 k}^{\Gamma}, X)$ lie on $[a_0, \infty)$. Moreover, if $m$ denotes the number of distinct zeros with odd multiplicity on $[a_0, \infty)$, then $m + c(\Gamma, \mathcal{F}) - s_{2 k}^1(\Gamma) \geqslant \deg(P(E_{2 k}^{\Gamma}, X))$.}\vspace{0.15in}

\section{Proof of Theorem \ref{th-main}}\label{sec-pr-th}
Now, we state the following theorem:
\begin{theorem}
Let $\Gamma$ be a genus zero group that is good for the weight $2 k$. Suppose that $\mathcal{F}$ is acceptable for $\Gamma$ and that $j_{\Gamma}$ has real Fourier coefficients. Then all but at most $c(\Gamma, \mathcal{F}) - s_{2 k}^1(\Gamma)$ of the zeros of $P(E_{2 k}^{\Gamma}, X)$ lie on $(- \infty, a_1]$. Moreover, if $m$ denotes the number of distinct zeros with odd multiplicity on $(- \infty, a_1]$, then $m + c(\Gamma, \mathcal{F}) - s_{2 k}^1(\Gamma) \geqslant \deg(P(E_{2 k}^{\Gamma}, X))$.\label{th-main1}
\end{theorem}

Then, by Theorem $\ref{th-main0}'$ and Theorem \ref{th-main1}, Theorem \ref{th-main} follows. The proof of Theorem \ref{th-main1} is similar to that of Theorem \ref{th-main0}, thus we will give only the outline of the proof:

\begin{proof}[Outline Proof of Theorem $\ref{th-main1}$]
We write $A = \coprod_{i=1}^{n} A_i$, where $A_i$ are smooth paths such that $j_{\Gamma} |_{\overline{A_i}}$. Similar to $\mathcal{F}_i$, we define closures $\mathcal{F}_i^{-}$ which satisfy $\overline{\mathcal{F}} = \bigcup_{i=1}^{n+2} \mathcal{F}_i^{-}$ and the properties for $\mathcal{F}_i$: $(1)$, $(2)$, $(3)$, and $(5)$ and
\begin{trivlist}
\item[$(4)'$] $\mathcal{F}_{n+1}^{-} = \overline{\mathcal{F}} \cup \{ z \, : \, \Re(z) \leqslant 0 \} \cup \{ z \, : \, \Im(z) \geqslant y_1 \}$ and $\mathcal{F}_{n+2}^{-} = \overline{\mathcal{F}} \cup \{ z \, : \, \Re(z) \geqslant 0 \} \cup \{ z \, : \, \Im(z) \geqslant y_1 \}$, where $y_1$ is defined as in definition $(\ref{y1})$
\end{trivlist}

For $z = x + i \, y \in \mathcal{F}_i^{-}$, define a path
\begin{equation}
\gamma_i^{-}(z) := \{ j_{\Gamma}^{-1}(\xi + i \, \Im(j_{\Gamma}(z))) \, : \, \xi \in (-\infty, \Re(j_{\Gamma}(z))] \},
\end{equation}
which is also traversed from $\infty$ to $z$, and let
\begin{equation}
\Phi_i^{-}(z) := \int_{\gamma_i^{-}(z)} \Re \left( P(f, j_{\Gamma_(s)}) \overline{P(g, j_{\Gamma_(s)})} |\Upsilon(s)|^2 t^{2k-2} \right) d t
\end{equation}
be the path integral where $\Im(s) = t$. Then, we obtain the following result: (cf. \cite[Proposition 4.1]{H})
\begin{proposition}
Let $\Gamma$ be good for the weight $2 k$ of genus zero and $\mathcal{F}$ be an acceptable fundamental domain for $\Gamma$. Let $f$ and $g$ be modular forms in $M_{2 k}(\Gamma)$ such that $f \, g$ vanishes at every cusp of $\Gamma$. Defining $\mathcal{F}_i^{-}$ and $\Phi_i^{-}$ as above, we have that
\begin{equation}
\iint_{\mathcal{F}} \Re \left( P(f, j_{\Gamma_(z)}) \overline{P(g, j_{\Gamma_(z)})} |\Upsilon(z)|^2 y^{2k-2} \right) d x d y
 = - \sum_{i = 1}^{n} \int_{A_i} \Phi_i^{-}(z) d x,
\end{equation}
where $z := x + i \, y$, and the $A_i$ are smooth paths traversed in an anti-clockwise direction.\label{prop-4.1-}
\end{proposition}
The proof of the above proposition is similar to that of Proposition 4.1 of \cite{H}. Note that both $\gamma_i$ and $\gamma_i^{-}$ are paths traversed from $\infty$ to $z$.

Let $\mathfrak{a}_1^{-}$, $\ldots$ , $\mathfrak{a}_{m'}^{-}$ be the zeros of $P(E_{2 k}^{\Gamma}, X)$ on $(-\infty, a_1]$ excluding images of cusps by $j_{\Gamma}$ that have odd multiplicity, and let
\begin{equation}
Q_1^{-}(X) := \prod_{i = 1}^{m'} (X - \mathfrak{a}_i^{-}) \prod_{j = 1}^{c(\Gamma, \mathcal{F}) - s_{2 k}^1(\Gamma, \mathcal{F})} (X - j_{\Gamma}(\mathfrak{b}_j)).
\end{equation}
We have $\langle E_{2 k}^{\Gamma}, Q_1^{-}(j_{\Gamma}) \Upsilon \rangle = 0$, and by Proposition \ref{prop-4.1-},
\begin{equation}
\Re (\langle E_{2 k}^{\Gamma}, Q_1^{-}(j_{\Gamma}) \Upsilon \rangle)
 = - \sum_{i = 1}^{n} \int_{A_i} \left( \int_{-\infty}^{j_{\Gamma}(z)} P(E_{2 k}^{\Gamma}, j_{\Gamma}) \, Q_1^{-}(j_{\Gamma}) \, (|\Upsilon|^2 t^{2 k}) (j_{\Gamma}) \, t^{-2} (j_{\Gamma}) \, t' (j_{\Gamma}) \, d j_{\Gamma} \right) d x.
\end{equation}
Since $P(E_{2 k}^{\Gamma}, X) Q_1^{-}(X) t'(X)$ is nonnegative or nonpositive on $(-\infty, j_{\Gamma}(z))$, we can then prove $m + c(\Gamma, \mathcal{F}) - s_{2 k}^1(\Gamma) \geqslant \deg(P(E_{2 k}^{\Gamma}, X))$.
\end{proof}

\quad\\
\appendix
\begin{center}
{\bfseries APPENDICES.}
\end{center}\quad

In Appendix \ref{sec-another-pf}, we demonstrate another proof of Theorem \ref{th-main1} using conjugate groups. In Appendix \ref{sec-no-acc}, we observe the location of the zeros of the Eisenstein series for some good groups with no acceptable fundamental domain.\\

\section{Another proof of Theorem \ref{th-main1}}\label{sec-another-pf}
Let $\Gamma$ be a Fuchsian group of the first kind with $\infty$ as a cusp, and let $\Gamma' = T_{h/2}^{-1} \, \Gamma \, T_{h/2}$ where $T_{h/2} = \left(\begin{smallmatrix} 1 & h/2 \\ 0 & 1 \end{smallmatrix}\right)$. For example, when $\Gamma = \Gamma_0(2)$, then we have $\Gamma' = \Gamma_0^{*}(4)$.

We can easily show that the map
\begin{equation*}
M_k(\Gamma) \ni f(z) \mapsto f(z + h/2) \in M_k(\Gamma')
\end{equation*}
is an isomorphism. Furthermore, we have
\begin{equation}
E_{2 k}^{\Gamma'}(z) = E_{2 k}^{\Gamma}(z + h/2). \label{EkGp}
\end{equation}
If $f \in M_{2 k}^{\infty}(\Gamma)$ has a Fourier expansion at infinity given by
\begin{equation*}
f(z) = \sum_{n=n_0}^{\infty} a_n q_h^n, \quad q_h := e^{2 \pi i  z / h}, \quad n_0 \in \mathbb{Z},
\end{equation*}
then since $q_h |_{z = z + h/2} = - q_h$ we have
\begin{equation*}
f(z + h/2) = \sum_{n=n_0}^{\infty} ((-1)^n a_n) q_h^n \; \in M_{2 k}^{\infty}(\Gamma').
\end{equation*}
\begin{remark}
For the weight $2 k$, if the group $\Gamma$ is good, then the group $\Gamma' = T_{h/2}^{-1} \, \Gamma \, T_{h/2}$ is also good.
\end{remark}
We also have $s_{2 k}^1(\Gamma') = s_{2 k}^1(\Gamma)$.

In addition, if canonical hauptmodul $j_{\Gamma}$ has a Fourier expansion of the form
\begin{equation*}
j_{\Gamma}(z) = \frac{1}{q_h} + \sum_{n=1}^{\infty} a_n q_h^n \; \in M_0^{\infty}(\Gamma),
\end{equation*}
then  
\begin{equation}
j_{\Gamma'}(z) = \frac{1}{q_h} + \sum_{n=1}^{\infty} ((-1)^{n-1} a_n) q_h^n = - j_{\Gamma}(z \pm h/2) \; \in M_0^{\infty}(\Gamma') \label{jGp}
\end{equation}
is a canonical hauptmodul for $\Gamma'$.

Let $\mathcal{F}$ be a fundamental domain for $\Gamma$ contained in $\{ z \, : \, -h/2 \leqslant \Re(z) < h/2 \}$, and write
\begin{align*}
\mathcal{F}^{-} &:= \mathcal{F} \cap \{ z \, : \, -h/2 \leqslant \Re(z) < 0 \},\\
\mathcal{F}^{+} &:= \mathcal{F} \cap \{ z \, : \, 0 \leqslant \Re(z) < h/2 \}.
\end{align*}
By the correspondence
\begin{equation}
\begin{split}
\mathcal{F}^{-} \; \ni \; z \; \mapsto \; z+h/2 \; \in \; \mathcal{F}^{-} + h/2,\\
\mathcal{F}^{+} \; \ni \; z \; \mapsto \; z-h/2 \; \in \; \mathcal{F}^{+} - h/2,
\end{split}\label{corr-1/2}
\end{equation}
we have that
\begin{equation}
\mathcal{F}' := (\mathcal{F}^{+} - h/2) \cup (\mathcal{F}^{-} + h/2)
\end{equation}
is a fundamental domain for $\Gamma'$, and we denote by $A'$ the lower arcs of $\mathcal{F}'$. By the above correspondence, we have the following fact.
\begin{remark}
If a fundamental domain $\mathcal{F}$ for $\Gamma$ is acceptable, then the domain $\mathcal{F}'$ in the above definition is an acceptable fundamental domain for $\Gamma'$. Furthermore, we have $c(\Gamma, \mathcal{F}) = c(\Gamma', \mathcal{F}')$.
\end{remark}

Now, let ${y_0}' := \inf\{ y \, : \, \pm h/2 + i \, y \in \partial \mathcal{F}' \} \; (= y_1)$ and ${a_0}' = \j_{\Gamma'}({y_0}') \; (= - a_1)$. Similarly to Theorem $\ref{th-main0}'$, we can show the following proposition:
\begin{proposition}
For the group $\Gamma'$ and the domain $\mathcal{F}'$ in the above definitions, all but at most $c(\Gamma', \mathcal{F}') - s_{2 k}^1(\Gamma')$ of the zeros of $P(E_{2 k}^{\Gamma'}, X)$ lie on $[{a_0}', \infty)$. Moreover, if $m$ denotes the number of distinct zeros with odd multiplicity on $[{a_0}', \infty)$, then $m + c(\Gamma', \mathcal{F}') - s_{2 k}^1(\Gamma') \geqslant \deg(P(E_{2 k}^{\Gamma'}, X))$.\label{prop-main1p}
\end{proposition}

By relations (\ref{EkGp}) and (\ref{jGp}), we have $P(E_{2 k}^{\Gamma'}, X) = (-1)^d P(E_{2 k}^{\Gamma}, -X)$ where $d = \deg(P(E_{2 k}^{\Gamma}, X))$. Thus, the above proposition is equivalent to Theorem \ref{th-main1}.\\

\section{On zeros of Eisenstein series for some good groups with no acceptable fundamental domain}\label{sec-no-acc}
Let $\Gamma$ be a Fuchsian group of the first kind with $\infty$ as a cusp. If $\Gamma$ does not have any acceptable fundamental domain, then the zeros of $E_{2 k}^{\Gamma}$ do not always lie on $\{ z \in \mathfrak{H} \, : \, j_{\Gamma}(z) \in \mathbb{R}\}$.

In this section, we consider the normalizers of the congruence subgroup $\Gamma_0(N)$. (see \cite{CN} and \cite{S}) For the levels $N \leqslant 12$, the normalizers with no acceptable fundamental domain are the following:
\begin{equation*}
\Gamma_0(5), \; \Gamma_0(6)+2, \; \Gamma_0(7), \; \Gamma_0(9), \; \Gamma_0(10)+2, \; \Gamma_0(10), \; \Gamma_0^{*}(11), \; \text{and} \; \Gamma_0(12)+3.
\end{equation*}

For $\Gamma_0(5)$, $\Gamma_0(6)+2$, $\Gamma_0(7)$, $\Gamma_0(10)+2$, $\Gamma_0(10)$, and $\Gamma_0^{*}(11)$, we can observe from numerical calculations that many of the zeros of the Eisenstein series for the $\infty$ cusp do not lie on the lower arcs of their fundamental domains. However, when the weight of Eisenstein series increases, then the location of the zeros seems to approach these lower arcs.

\begin{figure}[hbtp]
\begin{center}
{{\small $\mathcal{F}_{6+2}$}\includegraphics[width=1.5in]{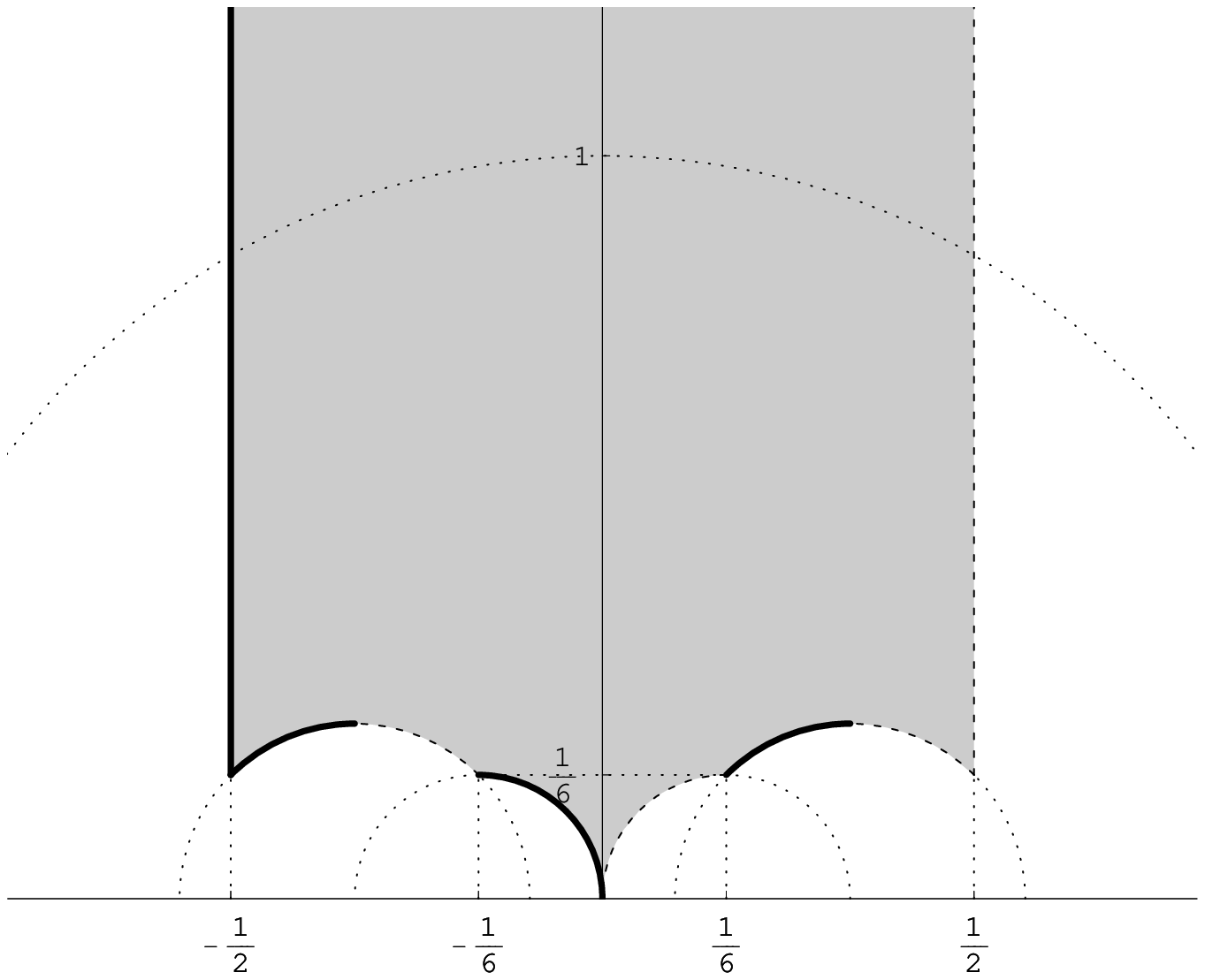}}\quad
{{\small Image of $A_{6+2}$ by $J_{6+2}$}\includegraphics[width=2in]{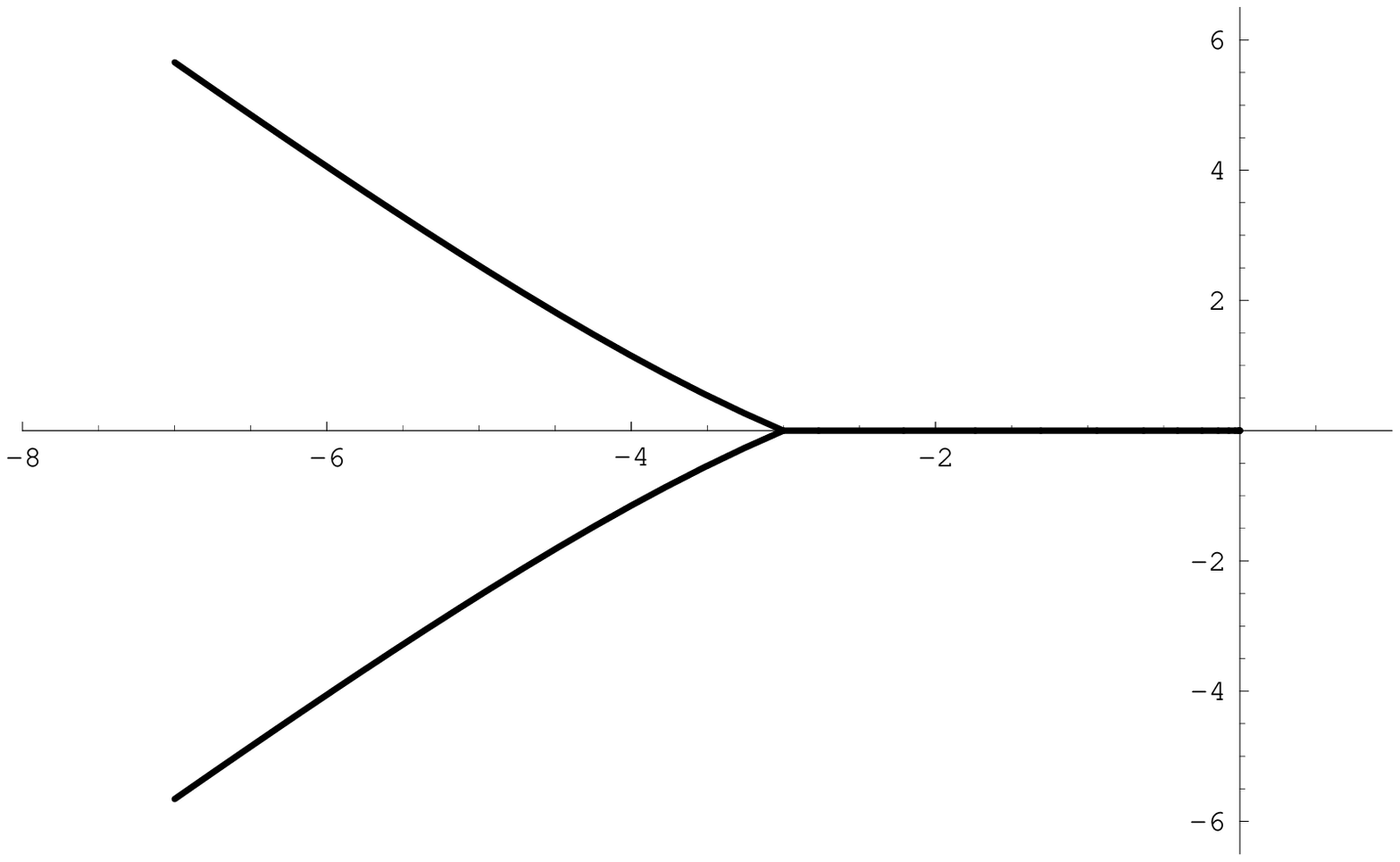}}\\
\quad\\
\includegraphics[width=2.5in]{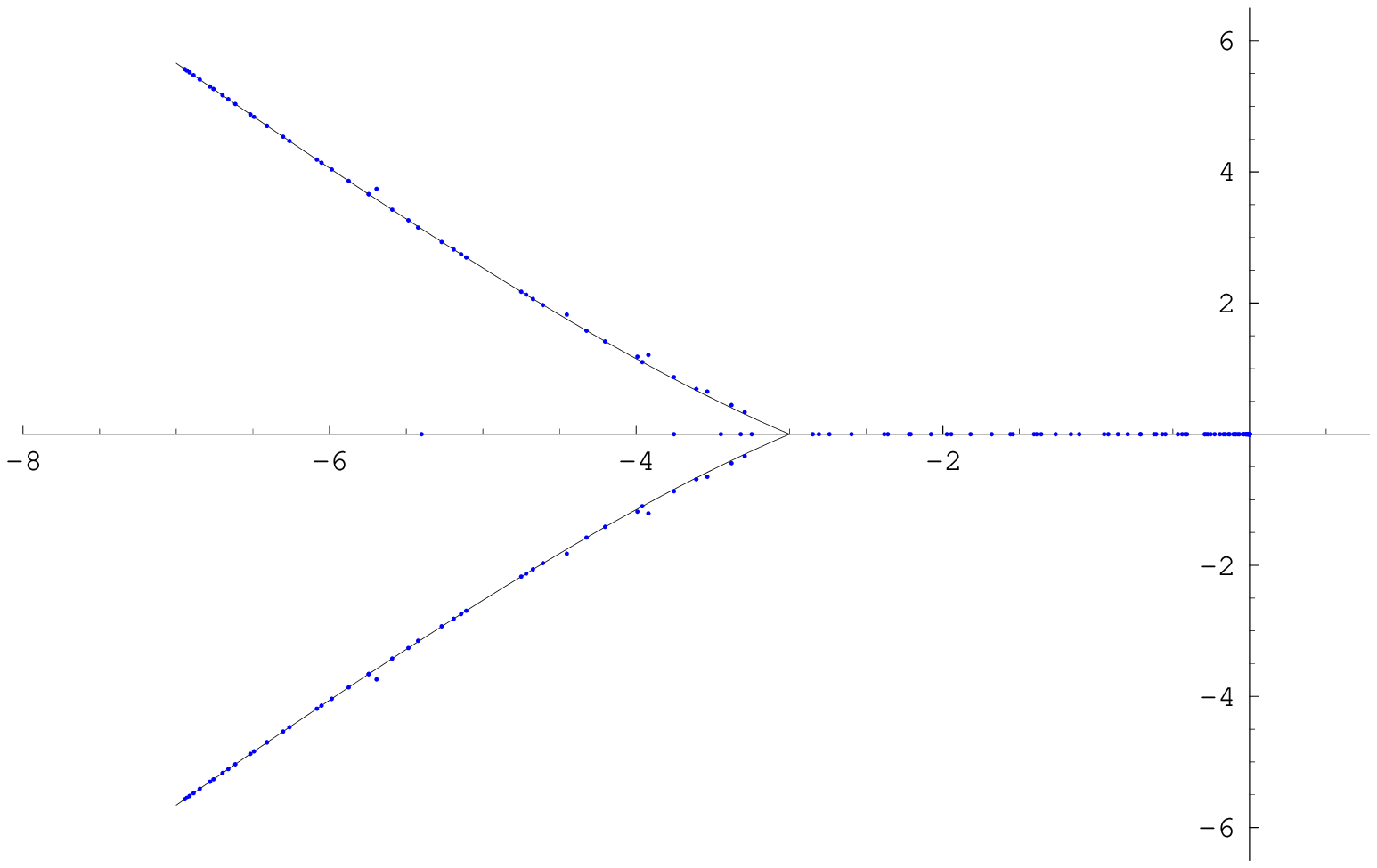}\\
{\small Image of the zeros of $E_{k}^{\Gamma_0(6)+2}$ for $4 \leqslant k \leqslant 40$ by $J_{6+2}$}\\
\quad\\
\includegraphics[width=2.5in]{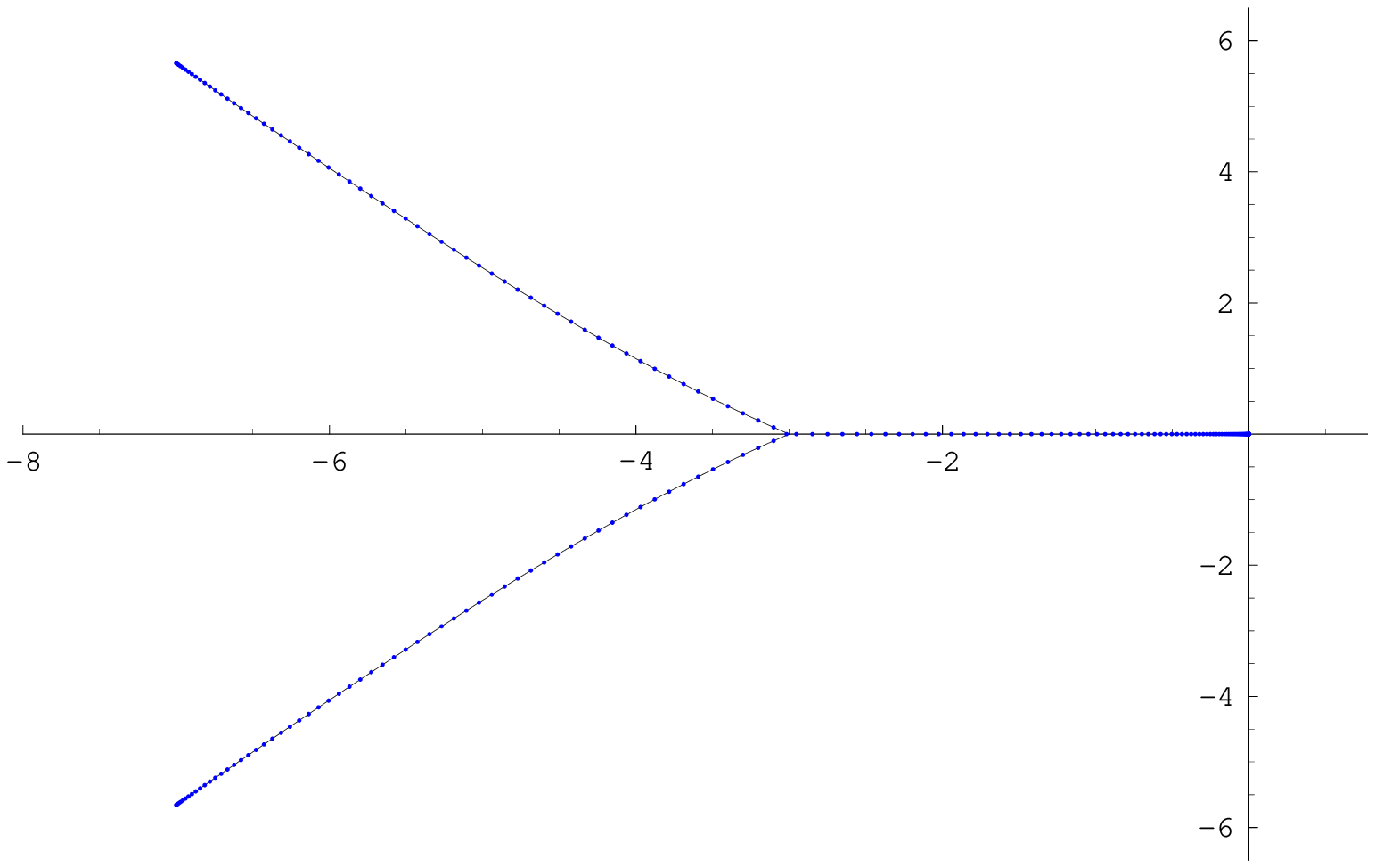}\\
{\small Image of the zeros of $E_{1000}^{\Gamma_0(6)+2}$ by $J_{6+2}$}
\end{center}
\caption{$\Gamma_0(6)+2$}\label{6D}
\end{figure}

In Figure \ref{6D}, we show graphs for $\Gamma_0(6)+2 = \Gamma_0(6) \cup W_{6, 2} \Gamma_0(6)$, where $W_{6, 2} := \left(\begin{smallmatrix} -\sqrt{2} & -1/\sqrt{2}\\ 3\sqrt{2} & \sqrt{2}\end{smallmatrix}\right)$. We denote by $\mathcal{F}_{6+2}$ its fundamental domain, by $A_{6+2}$ its lower arcs, and by $j_{6+2}$ its hauptmodul.

\begin{figure}[hbtp]
\begin{center}
{{\small $\mathcal{F}_9$}\includegraphics[width=1.5in]{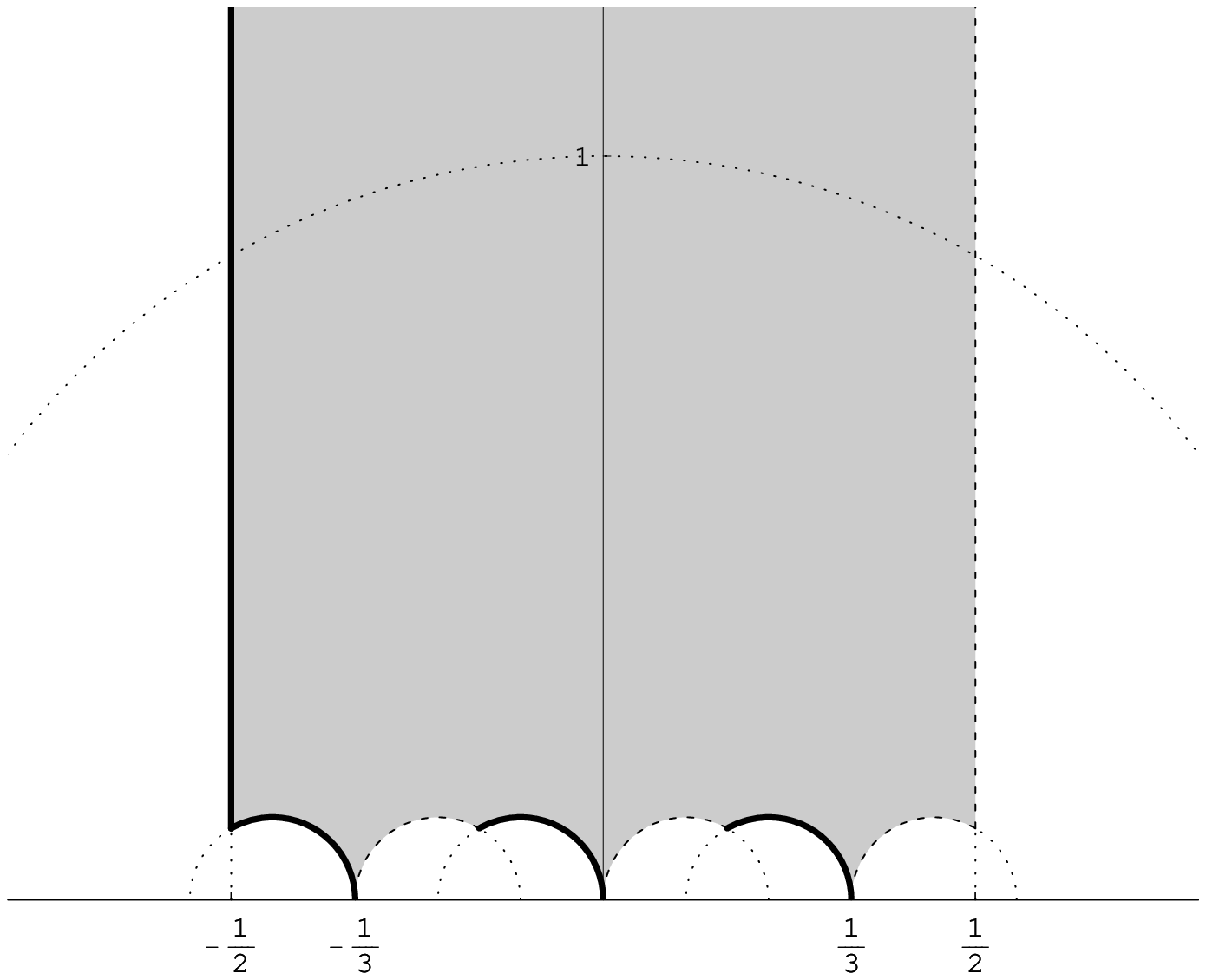}}\qquad
{{\small Image of $A_9$ by $J_9$}\includegraphics[width=1.5in]{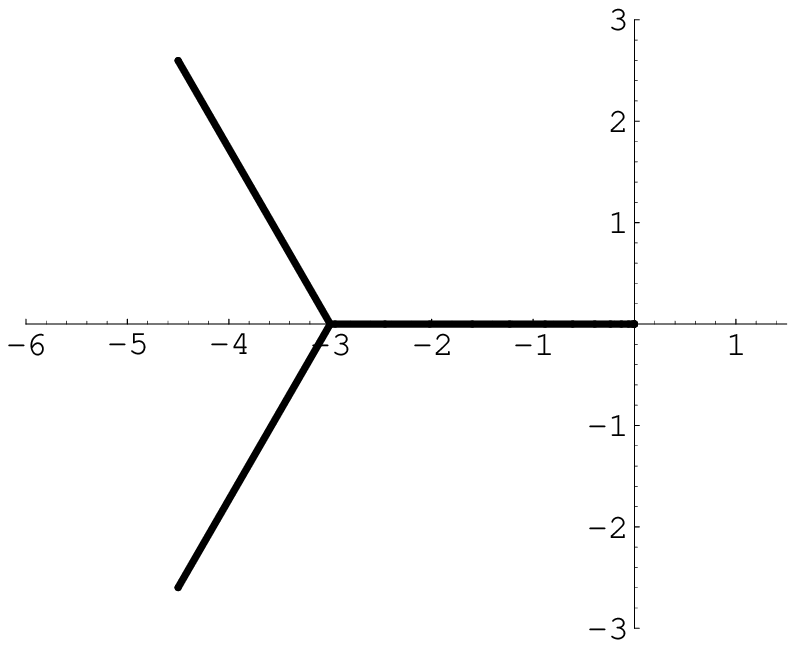}}
\end{center}
\caption{$\Gamma_0(9)$}\label{9B}
\end{figure}

On the other hand, $\Gamma_0(9)$ and $\Gamma_0(12)+3$ appear to be the special cases. For $\Gamma_0(9)$, we have
\begin{equation}
E_{2 k}^{\Gamma_0(9)}(z) = E_{2 k}^{\Gamma_0(3)}(3 z).
\end{equation}
Since by Theorem \ref{th-main} all of the zeros of $E_{2 k}^{\Gamma_0(3)}$ in its fundamental domain lie on its lower arcs, all of the zeros of $E_{2 k}^{\Gamma_0(9)}$ in its fundamental domain (see Figure \ref{9B}) lie on its lower arcs. We denote the fundamental domain of $E_{2 k}^{\Gamma_0(9)}$ by $\mathcal{F}_9$, its lower arcs by $A_9$, and its hauptmodul by $j_9$.

For $\Gamma_0(12)+3$, we have $E_{2 k}^{\Gamma_0(12)+3}(z) = E_{2 k}^{\Gamma_0(6)+3}(2 z)$. Thus, we can prove by Theorem \ref{th-main} that all but four of the zeros of $E_{2 k}^{\Gamma_0(12)+3}$ in its fundamental domain lie on its lower arcs. However, both $\Gamma_0(9)$ and $\Gamma_0(12)+3$ do not have any acceptable fundamental domains.

\quad\\

\end{document}